\documentclass[11pt]{article}
\usepackage{amssymb,amsfonts,amsthm}
\usepackage{color}
\usepackage[latin1]{inputenc}

\def\openC{{\rm C\kern-.18cm\vrule width.8pt height 7pt depth-.2pt \kern.18cm}}
\def\openN{{{\rm I}\kern-.16em {\rm N}}}
\def\openR{{{\rm I}\kern-.16em {\rm R}}}
\def\openT{{{\rm T}\kern-.42em {\rm T}}}
\def\openZ{{{\rm Z}\kern-.28em{\rm Z}}}

\def\la{\langle}
\def\ra{\rangle}

\def\eop{\hfill\rule{2.5mm}{2.5mm}}
\def\pf{\par\smallbreak\noindent {\bf Proof.} \ }

\newtheorem{thm}{Theorem}[section]
\newtheorem{cor}[thm]{Corollary}
\newtheorem{lem}[thm]{Lemma}
\newtheorem{prop}[thm]{Proposition}

\theoremstyle{definition}

\def\eop{\hfill\rule{2.5mm}{2.5mm}}
\begin{document}

\title{
{\textbf{\Large{Estimates for Fourier sums and eigenvalues of integral operators via multipliers on the sphere}}} \vspace{-4pt}
\author{
T. Jord\~{a}o\,\thanks{
Partially supported by FAPESP, grant $\#$ 2011/21300-7}
\,\,\&\,\,
V. A. Menegatto}
}
\date{}
\maketitle \vspace{-30pt}
\bigskip

\begin{center}
\parbox{13 cm}{{\small We provide estimates for weighted Fourier sums of integrable functions defined on the sphere when the weights originate from a multiplier operator acting on the space where the function belongs.\ That implies refined estimates for weighted Fourier sums
of integrable kernels on the sphere that satisfy an abstract H\"{o}lder condition based on a parameterized family of multiplier operators defining an approximate identity.\ This general estimation approach includes an important class of multipliers operators, namely, that defined by convolutions with zonal measures.\ The estimates are used to obtain decay rates for the eigenvalues of positive integral operators on $L^2(S^m)$ and generated by a kernel satisfying the H\"{o}lder condition based on multiplier operators on $L^2(S^m)$.}}
\end{center}


\thispagestyle{empty}

%
%

\section{Introduction}\label{s1}

In the recent paper \cite{jordao}, the use of a modulus of smoothness defined by the shifting operator to estimate Fourier coefficients of functions on the unit sphere $S^m$ in $\mathbb{R}^{m+1}$ has turned out to be extremely efficient in the deduction of decay rates for the sequence of eigenvalues of certain integral operators acting on spaces of integrable functions on $S^m$.\ Indeed, through a minor generalization of estimates originally obtained in \cite{ditzian}, one of the main results in \cite{jordao} deduces decay rates for the sequence of eigenvalues of integral operators generated by a Mercer's kernel satisfying a H\"{o}lder condition based on the shifting operator as introduced in \cite{samko}.\ Within the spherical setting, this result is an improvement upon classical results of the same type deduced in the early eighties in \cite{kuhn}.\ The ultimate target in the present paper is to obtain a result in this same framework, but using a H\"{o}lder condition defined by a parameterized family of bounded multiplier operators.

For the purpose of a formal presentation of the results, we need to introduce notation.\ First of all, we endow $S^{m}$ with its surface measure $\sigma_m$ and fix $m\geq 2$.\ For $p\geq 1$, we denote by $L^p(S^{m}):=L^p(S^{m},\sigma_m)$ the usual Banach space of integrable functions equipped with its $p$-norm $\|\cdot\|_p$ given by
$$
\|f\|_p=\left(\frac{1}{\omega_m}\int_{S^{m}}|f(x)|^pd\sigma_m(x)\right)^{1/p}, \quad f\in L^p(S^m),
$$
where $\omega_m$ is the surface area of $S^{m}$:
$$
\omega_m=\frac{2\pi^{(m+1)/2}}{\Gamma((m+1)/2)}.
$$

If we write $\la \cdot, \cdot \ra_2$ to denote the usual inner product of $L^2(S^{m})$, the orthogonal projection of $L^2(S^{m})$ onto the space $\mathcal{H}_k^m$ of all spherical harmonics of degree $k$ in $m+1$ variables will be written as ${\cal Y}_k$.\
If $\{Y_{k,j}: j=1,2,\ldots, d_k^m\}$ is an orthonormal basis of $\mathcal{H}_k^m$ with respect to $\la \cdot, \cdot \ra_2$, then
$$
{\cal Y}_k (f)(x)=\sum_{j=1}^{d_k^m}\hat{f}(k,j)Y_{k,j}, \quad x \in S^{m}.
$$
The Fourier coefficients in the above sum are computed through the formula
$$
\hat{f}(k,j)=\frac{1}{\omega_m} \int_{S^{m}} f(y)\overline{Y_{k,j}(y)}\, d\sigma_m(y), \quad j=1,2,\ldots, d_k^m.
$$
Additional information on the projections ${\cal Y}_k$ can be found in \cite{daixu} and references quoted there.

Next, we introduce multiplier operators of spherical harmonic expansions.\ A linear operator $T$ on $L^p(S^{m})$ is called a {\em multiplier operator} if there exists a sequence $\{\eta_k\}$ of complex numbers such that
\begin{equation} \label{alter}
{\cal Y}_k(Tf)=\eta_k {\cal Y}_k (f), \quad f \in L^p(S^{m}), \quad k=0,1,\ldots.
\end{equation}
The sequence $\{\eta_k\}$ is called the {\em sequence of multipliers} of $T$.\ An important category of multiplier operators are those which are invariant under rotations of $\mathbb{R}^{m+1}$ leaving a pole fixed, a typical example being the convolution with a zonal measure as exploited in \cite{dunkl}.\ Dunkl characterized the bounded multiplier operators on $L^1(S^m)$ as those given by such convolutions and, as far as we know, there is no similar characterization for multiplier operators on $L^2(S^m)$.\ However, it is not hard to see that a multiplier operator on $L^2(S^m)$ is bounded if and only if its sequence of multipliers is bounded.\ On the other hand, a bounded multiplier operator on $L^2(S^m)$ is self-adjoint and invariant under the group of rotations of $\mathbb{R}^{m+1}$ (\cite{daixu}).\ Classical results about inclusions among spaces of multipliers operators show that the class of multipliers operators on $L^2(S^m)$ is strictly bigger than the class of multiplier operators on $L^1(S^m)$.\ References on this particular subject are \cite{figa-talamanca, price, rieffel}.

A series of papers authored by Z. Ditzian and collaborators culminated with useful estimates for certain Fourier sums of functions in $L^p(S^m)$ $(1\leq p\leq 2)$ via a modulus of smoothness defined by certain combinations of the shifting operator on the sphere (\cite{ditzian}).\ The proof of those estimates required the following Hausdorff-Young type formulas ($q=$ the conjugate exponent of $p$):
$$
\left\{\sum_{k=0}^{\infty}(d_k^m)^{(2-q)/2q}\left[\sum_{j=1}^{d_k^m}|\hat{f}(k,j)|^2 \right]^{q/2}\right\}^{1/q} \leq \omega_{m}^{(p-2)/2p}\|f\|_p,\quad f \in L^p(S^m),\quad 1< p \leq 2,
$$
and
$$
\sup_{k\geq 0}\left\{(d_k^m)^{-1/2}\left[\sum_{j=1}^{d_k^m}|\hat{f}(k,j)|^2 \right]^{1/2}\right\}\leq \omega_{m}^{-1/2}\|f\|_1, \quad f \in L^1(S^m).
$$

In Section 2 of our paper, we provide similar estimates for certain sums of Fourier coefficients of an integrable function on $S^m$, replacing $f$ with $M f - f$, in which $M$ is a multiplier operator on $L^p(S^m)$.\ The coefficients of the internal sum in our inequalities depend not only on the dimensions $d_k^m$ but also on the pertinent sequence of multipliers of $M$.\ If the intention is to control the growth of the Fourier coefficients as $k \to \infty$, then one may think in the replacement of $M$ with a family $\{M_t: t \in (0,\pi)\}$ of multiplier operators which is an approximate identity in the sense that
$$
\lim_{t \to 0}\|M_t f -f\|_p =0,\quad f \in L^p(S^m).
$$
Keeping this in mind and restricting ourselves to the case $p=2$, we extend the estimates to
general kernels and also kernels that satisfy a H\"{o}lder assumption defined by a parameterized family of multiplier operators.\ Many commonly-used average operators can be used in the construction of a parameterized family of multiplier operators with the approximate identity property.\ In particular, the combinations of the shifting operator in \cite{dai} provide a solid example in the framework we consider.

In Section 3, we apply the results obtained in Section 2 in the deduction of decay rates for the sequence of eigenvalues of integral operators on the sphere, in the case the operator is generated by a Mercer-like kernel satisfying an abstract H\"{o}lder condition defined by a parameterized family of multipliers operators on $L^2(S^m)$.\ As mentioned before, a similar technique was firstly used in \cite{jordao} to deduce decay rates for similar integral operators, those having the generating kernel satisfying a H\"{o}lder condition defined by the shifting operator.\ The decay obtained in that setting coincided with that one deduced by K\"{u}hn in \cite{kuhn} in the case the generating kernel satisfies a general integrated H\"{o}lder condition (this condition encompasses the standard one).\ In addition to the method we use, a key contribution in our paper resides in the use of an abstract H\"{o}lder condition that can be applied in many other situations, including those in \cite{jordao,kuhn}, at least when one considers either the spherical setting or a similar one.\ Despite getting the very same decay rates, our arguments do not require the continuity of $K$ explicitly.\ On the other hand, our setting includes Mercer-like kernels defined by a parameterized family in which all elements are convolutions with zonal measures, a category that includes many operators defined by an average process.\

Section 4 contains concrete examples that illustrate our findings.

\section{Estimates of Fourier-like sums}

This section begins with the deduction of basic inequalities for sums of Fourier coefficients of integrable functions and kernels, generalizations of both Theorem 6.1 in \cite{ditzian} and the Hausdorff-Young inequalities mentioned in the previous section.\ If one replaces the sphere with the Euclidean space where it sits and the Fourier sums with the standard Fourier transformation, then the results are comparable to those proved in \cite{bray}.\ A differential in our favor is the fact that we do not need to make use of either $K$-functionals or moduli of smoothness in the arguments leading to the inequalities.\ At the end of the section, we estimate upon certain Fourier sums of kernels satisfying an abstract H\"{o}lder defined by a parameterized family of multiplier operators.

\begin{thm} \label{th1} Let $M$ be a multiplier operator on $L^p(S^{m})$ with corresponding sequence of multipliers $\{\eta_k\}$.\ If $p \in (1,2]$, then
$$
\left\{\sum_{k=1}^{\infty}(d_k^m)^{(2-q)/2q}|\eta_k-1|^{q} \left[\sum_{j=1}^{d_k^m}|\hat{f}(k,j)|^2\right]^{q/2}\right\}^{1/q}\leq \omega_{m}^{(p-2)/2p}\|M f-f\|_p,\quad f \in L^p(S^{m}),
$$
in which $q$ is the conjugate exponent of $p$.\ The inequality above becomes an equality in the case $p=2$.\ If $p=1$, then
$$
\sup_{k\geq 0}\left\{(d_k^m)^{-1/2}|\eta_k-1|\left[\sum_{j=1}^{d_k^m}|\hat{f}(k,j)|^2 \right]^{1/2}\right\}\leq \omega_{m}^{-1/2}\|M f-f\|_1,\quad f \in L^p(S^{m}).
$$
\end{thm}

\pf Fixing $f\in L^p(S^m)$, the linearity of the orthogonal projections and (\ref{alter}) imply that
$$
{\cal Y}_k(M f -f) =(\eta_k-1){\cal Y}_k(f),  \quad k\in \mathbb{Z}_+,
$$
whence
$$
\sum_{j=1}^{d_k^m}\widehat{(M f - f)}(k,j)Y_{k,j}=(\eta_k-1)\sum_{j=1}^{d_k^m}\hat{f}(k,j) Y_{k,j},\quad k\in \mathbb{Z}_+.
$$
Since the sums define polynomial functions, then we can compute the $L^2$-norm of both sides to obtain
$$
\sum_{j=1}^{d_k^m}\left|\widehat{(M f - f)}(k,j)\right|^2=|\eta_k-1|^2 \sum_{j=1}^{d_k^m}\left|\hat{f}(k,j)\right|^2,\quad k\in \mathbb{Z}_+,
$$
that is,
$$
(d_k^m)^{(2-q)/2q}\left(\sum_{j=1}^{d_k^m}\left|\widehat{(M f - f)}(k,j)\right|^2\right)^{q/2}=(d_k^m)^{(2-q)/2q}|\eta_k-1|^q \left(\sum_{j=1}^{d_k^m}\left|\hat{f}(k,j)\right|^2\right)^{q/2}.
$$
Applying the Hausdorff-Young formula we reach the first inequality in the statement of the theorem.\ As for the equality assumption in the case $p=2$, it suffices to apply Parseval's identity in the equality above.\ The inequality in the case $p=1$ is settled in a similar fashion.\eop

Next, we apply the equality in the previous theorem in order to obtain a similar result for square integrable kernels.\ If $K$ is a kernel in $L^2(S^{m}\times S^{m}):=L^2(S^{m}\times S^{m}, \sigma_m\times\sigma_m)$ with spherical harmonics expansion
\begin{equation}\label{kern}
K(x,y)=\sum_{k=0}^{\infty}\sum_{j=1}^{d_k^m}a_{k,j}Y_{k,j}(x)\overline{Y_{k,j}(y)}, \quad x,y \in S^{m},
\end{equation}
then every function $K^y$, $y \in S^m$, defined by
$$
K^y(x)=K(x,y), \quad x \in S^{m},
$$
belongs to $L^2(S^{m})$.\ In addition,
$$
\widehat{K^y}(k,j)=a_{k,j}\overline{Y_{k,j}(y)},\quad j=1,2,\ldots, d_k^m,\quad k\in \mathbb{Z}_+.
$$
Consequently,
$$
\frac{1}{\omega_m}\int_{S^{m}}\sum_{j=1}^{d_k}\left|\widehat{K^y}(k,j)\right|^2d\sigma(y)=\sum_{j=1}^{d_k^m}|a_{k,j}|^2,\quad k\in \mathbb{Z}_+.
$$
The result below is now evident since it can be obtained from this by integration of both sides of the equality in Theorem \ref{th1}.

\begin{thm} \label{multeq} Let $M$ be multiplier operator on $L^2(S^m)$ with multiplier sequence $\{\eta_k\}$.\ If $K$ is a kernel as in (\ref{kern})
then
$$
\sum_{k=0}^{\infty}|\eta_k-1|^2 \sum_{j=1}^{d_k^m}|a_{k,j}|^2=\frac{1}{\omega_m}\int_{S^m}\|M (K^y)-K^y\|_2^2\,d\sigma_m(y).
$$
\end{thm}

Next, we introduce a H\"{o}lder condition attached to a sequence $\{M_t: t \in (0,\pi)\}$ of multiplier operators on $L^2(S^m)$ and specialize the previous theorem to the case in which the kernel $K$ satisfies such a H\"{o}lder condition.\ We say that a kernel $K$ from $L^{2}(S^{m}\times S^{m})$ is {\em $\{M_t:t\in (0,\pi)\}$-H\"{o}lder} if there exist a real number $\beta \in (0,2]$ and a constant $B>0$ so that
\begin{equation}\label{intholder}
\int_{S^m}|M_t(K^y)(y)-K^y(y)|d\sigma_m(y) \leq B t^{\beta}.
\end{equation}
Clearly, this H\"{o}lder condition is implied by the more classical one which demands the existence of $\beta \in (0,2]$ and a
function $B$ in $L^1(S^{m})$ such that
\begin{equation}\label{SScondition}
\sup_x |M_t(K^y)(x)-K^y(x)| \leq B(y)t^{\beta},\quad y \in S^{m}, \quad t\in (0,\pi).
\end{equation}
The number $\beta$ appearing in the above definitions will be termed {\em the H\"{o}lder exponent of $K$ with respect to the family $\{M_t:t\in (0,\pi)\}$}.

At this point, we also need the notion of {\em $L^2$-positive definiteness}.\ For a kernel representable as in (\ref{kern}), it corresponds to
$$
a_{k,j} \geq 0, \quad j=1,2,\ldots, d_k^m, \quad k \in \mathbb{Z}_+.$$
In particular, the formula
$$K_{1/2}(x,y):=\sum_{k=0}^{\infty}\sum_{j=1}^{d_k^m}a_{k,j}^{1/2}Y_{k,j}(x)\overline{Y_{k,j}(y)}, \quad x,y \in S^{m},
$$
defines a positive definite element of $L^2(S^{m} \times S^{m})$ for which
\begin{equation}\label{sqroot}
\frac{1}{\omega_m}\int_{S^{m}} K_{1/2}(x,y)K_{1/2}(w,x) d\sigma_m(x) = K(w,y), \quad y,w \in S^{m}.
\end{equation}

In the sequel, the uniform boundedness of a family $\{M_t: t \in (0,\pi)\}$ of operators acting on $L^2(S^{m})$ will refer to the uniform boundedness of the numerical set $\{\|M_t\|: t\in (0,\pi)\}$.

\begin{lem} \label{lema} Let $\{M_t: t \in (0,\pi)\}$ be a uniformly bounded family of multiplier operators on $L^2(S^{m})$ with corresponding sequences of multipliers $\{\eta_k^t\}$.\ If $K$ is a $L^2$-positive definite $\{M_t:t\in (0,\pi)\}$-H\"{o}lder kernel, then there exists $C>0$ such that
$$\frac{1}{\omega_m}\int_{S^{m}}\|M_t(K_{1/2}^y)-K_{1/2}^y\|_2^2\,d\sigma(y)\leq C t^{\beta}, \quad t\in (0,\pi),$$
in which $\beta$ is the H\"{o}lder exponent of $K$ with respect to $\{M_t:t\in (0,\pi)\}$.
\end{lem}
\pf Assume $K$ has the features listed in the statement of the lemma.\ If $y \in S^{m-1}$ and $t \in (0,\pi)$, first observe that
\begin{eqnarray*}
\|M_t(K_{1/2}^y)-K_{1/2}^y\|_2^2 & = & \frac{1}{\omega_m} \int_{S^{m}}M_t(K_{1/2}^y)(z)\overline{M_t(K_{1/2}^y(z)}\, d\sigma_m(y)\\
     &  & \hspace*{15mm}-\frac{1}{\omega_m}\int_{S^{m}}M_t(K_{1/2}^y)(z)\overline{K_{1/2}^y(z)}\, d\sigma_m(y)\\
     &  & \hspace*{25mm}- \frac{1}{\omega_m}\int_{S^{m}}K_{1/2}^y(z)\overline{M_t(K_{1/2}^y)(z)}\, d\sigma_m(y)\\
     &  & \hspace*{40mm}+ \frac{1}{\omega_m}\int_{S^{m}}K_{1/2}^y(z)\overline{K_{1/2}^y(z)}\, d\sigma_m(y).
\end{eqnarray*}
The introduction of the Fourier expansions of the functions involved, some convenient calculations and the use of (\ref{sqroot}) lead to
$$
\|M_t(K_{1/2}^y)-K_{1/2}^y\|_2^2  =  M_t(M_t(K^y)-K^y)(y)-(M_t(K^y)-K^y)(y).
$$
Since each $M_t$ is self-adjoint, it is promptly seen that
\begin{eqnarray*}
\int_{S^{m}}M_t(M_t(K^y)-K^y)(y)d\sigma_m(y) & = & \int_{S^{m}}(M_t(K^y)-K^y)(y)M_t(1)(y)d\sigma(y)\\
                                             & = & \eta_0^t \int_{S^{m}}(M_t(K^y)-K^y)(y)d\sigma_m(y),\quad t \in (0,\pi).
                                             \end{eqnarray*}
It is now clear that
$$
\frac{1}{\omega_m} \int_{S^{m}}\|M_t(K_{1/2}^y)-K_{1/2}^y\|_2^2d\sigma_m(y)\leq  B(|\eta_0^t|+1)  t^{\beta}, \quad t \in (0,\pi).
$$
However,
$$|\eta_0^t| \leq \sup_k |\eta_k^t| \leq \|M_t\| \leq \sup\{\|M_s\|: s \in (0,\pi)\},\quad t\in(0,\pi),$$
and the inequality in the statement of the lemma follows.\eop

\begin{thm} \label{keyabst} Let $\{M_t: t \in (0,\pi)\}$ be a uniformly bounded family of multiplier operators on $L^2(S^{m})$, with corresponding multiplier sequences $\{\eta_k^t\}$.\ If $K$ is a $L^2$-positive definite $\{M_t:t\in (0,\pi)\}$-H\"{o}lder kernel, then there exists $C>0$ such that
$$\sum_{k=0}^{\infty}|\eta_k^t-1|^2 \sum_{j=1}^{d_k^m}a_{k,j} \leq C t^{\beta}, \quad t\in (0,\pi),$$
in which $\beta$ is the H\"{o}lder exponent of $K$ with respect to $\{M_t:t\in (0,\pi)\}$.
\end{thm}
\pf If $K$ is positive definite, then the same is true of $K_{1/2}$ and Theorem \ref{multeq} implies that
$$\sum_{k=0}^{\infty}|\eta_k^t-1|^2 \sum_{j=1}^{d_k^m}a_{k,j}=\frac{1}{\omega_m}\int_{S^m}\|M_t (K_{1/2}^y)-K_{1/2}^y\|_2^2\,d\sigma_m(y), \quad y\in S^{m} \quad t \in (0,\pi).$$
An application of Lemma \ref{lema} closes the proof.\eop

\section{Application to decay rates of eigenvalues}

Here, we employ the very last result in the previous section to extract decay rates for the sequence of eigenvalues of positive integral operators acting on $L^2(S^{m})$.\ By an integral operator we mean a bounded linear operator $\mathcal{K}:L^2(S^{m}) \to L^2(S^{m})$ having the form
$$\mathcal{K}(f)=\frac{1}{\omega_m}\int_{S^{m}}K(\cdot ,y)f(y)d\sigma_m(y),\quad f \in L^2(S^{m}),$$
in which $K$ is an element of $L^2(S^{m} \times S^{m})$ (the generating kernel of $K$).\ The integral operator $K$ is {\em positive}
whenever it is self-adjoint and the generating kernel $K$ is $L^2$-positive definite in the sense explained in the previous section.\ Basic spectral theory asserts that a positive integral operator $\mathcal{K}$ has at most countably many eigenvalues which can be ordered in a decreasing manner, say,
$$
\lambda_1(\mathcal{K})\geq \lambda_2(\mathcal{K})\geq \cdots \geq 0,
$$
with multiplicities being included.

If we use the expansion (\ref{kern}) for the generating kernel $K$ of a positive integral operator $\mathcal{K}$, then it is easily seen that
$$\mathcal{K}(Y_{k, j})=a_{k, j}, \quad j=1,2,\ldots, d_k^m, \quad k\in\mathbb{Z}_+.$$
In particular, the set $\{a_{k, j}: j=1,2,\ldots, d_k^m;k=0,1,\ldots\}$ is the set of eigenvalues of $\mathcal{K}$.\ Without any loss of generality, we can assume that for every $k$, $a_{k, 1}\geq a_{k, 2} \geq \cdots \geq a_{k, d_k^m}$.\ In particular, taking into account the ordering mentioned above,
$$\lambda_{d_n^{m+1}}(\mathcal{K})=\lambda_{1+d_1^m+\cdots +d_n^m}(\mathcal{K})=a_{nd_n^m},\quad n\in \mathbb{Z}_+.$$
On the other hand, if $K$ is to carry a smoothness assumption, then we can also assume that $a_{nj}\leq a_{kl}$ whenever $j=1,2,\ldots, d_k^m$, $l=1,2,\ldots, d_n^m$ and $n\geq k$ (\cite{seeley}).\ With all this in mind, we can now prove the general theorem below.

The theorem itself depends on sequences with special features in the way we now explain.\ A double indexed sequence $\{b_{k,n}\}$ of nonnegative real numbers is {\em half-bounded away from 0} if $\lim_{n\to \infty} b_{k,n}=0$, $k \in \mathbb{Z}_+$ and there exists a positive real number $M$ so that $b_{k,n} \geq M$, $k\geq n$.\ Roughly speaking, a half-bounded away from 0 double sequence can be ``controlled" in both cases, for small $k$ and large $k$, for all $n$.\ By looking at concrete cases (for example, a family of multiplier operators defined by convolutions), it is promptly seen that the requirement $\lim_{n\to \infty} b_{k,n}=0$, $k \in \mathbb{Z}_+$, in the above definition corresponds to the fact that the family is an approximate identity in $L^2(S^m)$ (see \cite{menega} for details).\ A simple example of a half-bounded away from 0 double sequence is $\{k/(k+n)\}$.

\begin{thm} \label{mainn} Let $\{M_t: t \in (0,\pi)\}$ be a uniformly bounded family of multiplier operators on $L^2(S^{m})$, with corresponding multiplier sequences $\{\eta_k^t\}$.\ Let $\mathcal{K}$ be an integral operator generated by a
$L^2$-positive definite $\{M_t:t\in (0,\pi)\}$-H\"{o}lder kernel $K$.\ If $\{|\eta_k^{1/n}-1| \}$ is half-bounded away from 0, then $\lambda_n(\mathcal{K})=O(n^{-1-\beta/m})$, as $n\to \infty$, in which $\beta$ is the H\"{o}lder exponent of $K$ with respect to $\{M_t:t\in (0,\pi)\}$.
\end{thm}
\pf Pick $C>0$ so that $|\eta_k^{1/n}-1| \geq C$, $k \geq n$.\ Since an application of Theorem \ref{keyabst} yields
$$\sum_{k=n}^{\infty}|\eta_k^{1/n}-1|^2 \sum_{j=1}^{d_k^m} a_{k, j}\leq C n^{-\beta}, \quad n=1,2,\ldots,$$
we can deduce that
$$\sum_{k=n}^{\infty} \sum_{j=1}^{d_k^m}a_{k, j}\leq C_1 n^{-\beta}, \quad n=1,2, \ldots,$$
for some $C_1>0$.\ Consequently,
$$d_n^m \sum_{k=n}^{\infty} a_{k, d_k^m} \leq \sum_{k=n}^{\infty} d_k^m a_{k, d_k^m}\leq C_1 n^{-\beta}, \quad n=1,2,\ldots.$$
Using the equivalence $d_n^m \asymp n^{m-1}$, as $n \to \infty$, we are reduced ourselves to an inequality of the form
$$\sum_{k=n}^{\infty} a_{k, d_k^m} \leq C_2 n^{-\beta-m+1}, \quad n=1,2,\ldots,$$
with $C_2>0$.\ Another estimation leads to
$$n^{\beta+m} a_{n, d_n^m} \leq n^{\beta+m-1}\sum_{k=n}^{\infty} a_{k, d_k^m} \leq C_2, \quad n=1,2,\ldots.$$
that is, $\lambda_{d_n^{m+1}}(\mathcal{K})=O(n^{-\beta-m})$, as $n\to \infty$.\ But, elementary calculations with this information implies that
$\lambda_n(\mathcal{K})=O(n^{-1-\beta/m})$ as $n \to \infty$.\eop

Several important multiplier operators, those given by certain averages on $S^m$ (see examples in the next section), present the following feature for the sequence of multipliers: there exists $s>0$ so that
$$|\eta_k^t-1|\asymp (\min\{1,tk\})^{s}, \quad t \in (0,\pi),\quad k \in \mathbb{Z}_+.$$
In other words, there exist positive constants $c_1$ and $c_2$ so that
$$c_1 (\min\{1,tk\})^{s} \leq |\eta_k^t-1| \leq c_2 (\min\{1,tk\})^{s}, \quad t \in (0,\pi),\quad k \in \mathbb{Z}_+.$$
Obviously, for such a sequence, the double sequence $\{|\eta_k^{1/n}-1|\}$ is half-bounded away from 0 and the following corollary holds.

\begin{cor} \label{mainncor} Let $\{M_t: t \in (0,\pi)\}$ be a uniformly bounded family of multiplier operators on $L^2(S^{m})$, with corresponding multiplier sequences $\{\eta_k^t\}$.\ Let $K$ be an integral operator generated by a positive definite $\{M_t:t\in (0,\pi)\}$-H\"{o}lder kernel $K$.\ If $$|\eta_k^{t}-1|  \asymp (\min\{1,tk\})^{s},\quad t \in (0,\pi), \quad k \in\mathbb{Z}_+,$$
for some $s>0$, then $\lambda_n(\mathcal{K})=O(n^{-1-\beta/m})$, as $n\to \infty$, in which $\beta$ is the H\"{o}lder exponent of $K$ with respect to $\{M_t:t\in (0,\pi)\}$.
\end{cor}

\section{A concrete case: convolutions with zonal measures}

This section contemplates several examples of families of multiplier operators that fit in the settings of Theorem \ref{mainn} and its corollary.\ The main class of multiplier operators we intend to consider is that of convolution operators with a family of zonal measures.

Let $SO_{m+1}$ be the group of rotations of $S^m$ and $G_x:= \{\mathcal{O}\in G: \mathcal{O}(x)=x\}$ the closed subgroup of $G$ that fixes a particular element $x$ of $S^m$.\ The set $\mathcal{M}(S^m)$ of all finite regular measures on $S^m$ becomes a Banach space when we define the norm of an element of $\mathcal{M}(S^m)$ as being its total variation.\ The set
$$
M_x(S^m):=\{\mu\in \mathcal{M}(S^m): \mu\circ \mathcal{O}=\mu, \,\,  \mathcal{O}\in G_x\}
$$
being a closed subspace of $\mathcal{M}(S^m)$, is likewise a Banach space under the same norm.\ If $x$ and $\varepsilon$ are elements of $S^m$ then $\mathcal{M}_x(S^m)$ and $\mathcal{M}_\varepsilon(S^m)$ are isomorphic.\ Indeed, if $\mathcal{O}_x^{\varepsilon} \in SO_{m+1}$ satisfies $\mathcal{O}_x^\varepsilon(x)=\varepsilon$ then the formula
$$
\varphi_x(\mu)= \mu\circ \mathcal{O}_x^\varepsilon, \quad \mu\in \mathcal{M}_p(S^m),
$$
defines a canonical isomorphism from $\mathcal{M}_\varepsilon(S^m)$ to $\mathcal{M}_x(S^m)$.\ In what follows, the total variation of an element $\mu$ of $\mathcal{M}(S^m)$ will be written as $|\mu|$.\ Hence, the norm of an element $\mu$ in either $\mathcal{M}(S^m)$ or $\mathcal{M}_x(S^m)$ is just
$$
|\mu|(S^m)=\sup\left\{\frac{1}{\omega_m}\left|\int_{S^m}fd\mu\right|: f\in L^1(S^m,\mu);\,|f|\leq 1\right\}.
$$
If $\mu$ is a positive element of $\mathcal{M}(S^m)$ then $|\mu|=\mu$ (\cite[p.85-87]{folland}).\ Finally, a procedure as above ratifies that
$$
L^2_x(S^m):=\{f\in L^2(S^m): f \circ \mathcal{O}=f,\,\,  \mathcal{O}\in G_x\}
$$
is a Banach space with the norm inherited from $L^2(S^m)$.

\begin{prop} \label{convo} (\cite{berens,dunkl})
Let $\varepsilon$ be a fixed pole in $S^m$.\ If $f$ belongs to $L^2(S^m)$ and $\mu$ is an element of $\mathcal{M}_\varepsilon(S^m)$, then the formula
\begin{equation}\label{convolution}
(f\ast\mu)(x):=\frac{1}{\omega_m}\int_{S^m}f(y)d\varphi_x(\mu)(y),
\end{equation}
defines an element of $L^2(S^m)$ satisfying $\|f\ast \mu\|_2\leq \|f\|_2|\mu|$.
\end{prop}

We will call $f\ast\mu$ {\em the spherical convolution of $f$ and $\mu$}.\ It is not hard to see that
$$
{\cal Y}_k(f\ast\mu)=\mu_k{\cal Y}_k(f), \quad f\in L^2(S^m),\quad \mu\in \mathcal{M}_\varepsilon(S^m),
$$
where
$$
\mu_k=\frac{1}{\omega_m}\int_{S^m}C_k^{(m-1)/2}(\varepsilon\cdot y)d\mu(y)
$$
and $C_k^{(m-1)/2}$ is the usual Gegenbauer polynomial associated with the real number $(m-1)/2$.\ In particular, the spherical convolution $f\ast\mu$ is, indeed, a multiplier operator.\ Additional information on the material described above can be found
in \cite{dunkl}.

Below is a list of examples in which the family of multiplier operators fittting into the main results of the previous section is defined by spherical convolutions with measures.
\vspace*{2mm}

\noindent \textbf{Shifting operator.} The usual \emph{shifting operator} is defined by the formula (\cite{berens,rudin})
$$
S_tf(x)=\frac{1}{R_m(t)}\int_{R_x^t}f(y)d\sigma_r(y),\quad x \in S^m, \quad f\in L^2(S^m),\quad t\in (0,\pi),
$$
in which $d\sigma_r(y)$ is the volume element of the rim $R_x^t:=\{y \in S^m: x\cdot y=\cos t\}$ and $R_m(t)=\omega_{m-1}(\sin t)^{m-1}$ is its total volume.\ If $A$ is a measurable subset of $S^m$, the formula
$$
\tilde{\mu}_t^m(A)=\omega_m\omega_{m-1}(A\cap R_\varepsilon^t),\quad t \in (0,\pi),
$$
defines an element $\mu_t^m:=R_m^{-1}(t)\tilde{\mu}_t^m$ of $\mathcal{M}_\varepsilon(S^m)$ for which
$$
S_t(f)=f\ast \mu_t^m, \quad f \in L^2(S^m), \quad t\in (0,\pi).
$$
Since
$$
{\cal Y}_k(S_t f)=\frac{C_k^{(m-1)/2}(\cos t)}{C_k^{(m-1)/2}(1)}{\cal Y}_k,\quad k \in \mathbb{Z}_+, \quad t \in (0,\pi),
$$
the sequence of multipliers of $S_t$ are obtained from
$$
\eta_k^t=\frac{C_k^{(m-1)/2}(\cos t)}{C_k^{(m-1)/2}(1)},\quad k \in \mathbb{Z}_+, \quad t \in (0,\pi).
$$
In Lemma 2.4 in \cite{belinsky}, it is proved that
$$ 0 < c_1 k^2 t^2 \leq 1-\frac{C_k^{(m-1)/2}(\cos t)}{C_k^{(m-1)/2}(1)} \leq c_2 k^2 t^2, \quad 0 < kt \leq \pi, \quad t\in (0,\pi/2],
$$
and that, for any $\tau >0$,
$$
\frac{C_k^{(m-1)/2}(\cos t)}{C_k^{(m-1)/2}(1)} \leq \alpha <1, \quad kt \geq \tau >0, \quad t\in (0,\pi/2],
$$
in which $c_1$ and $c_2$ are positive constants depending on $m$ and $\alpha$ is a constant depending on $m$ and $\tau$.\ As a consequence, it is seen that
$$
1-\frac{C_k^{(m-1)/2}(\cos t)}{C_k^{(m-1)/2}(1)} \asymp (\min\{1, kt\})^{2}, \quad k \in \mathbb{Z}_+, \quad t \in (0,\pi).
$$
The inequality
$$
\|S_t f\|_2 \leq \|f\|_2, \quad f \in L^2(S^m), \quad t \in (0,\pi),
$$
is all that is needed in order to see that $\{S_t: t \in (0,\pi)\}$ is uniformly bounded by 1.
\vspace*{2mm}

\noindent \textbf{Combinations of shiftings.} This example was developed in \cite{dai} and it can be seen an extension of the previous one.\ For $l=1,2,\ldots$, let $S_{l,t}$ be the operator on $L^2(S^m)$
 given by
$$
S_{l,t} (f)=-2{2l\choose l}^{-1}\sum_{j=1}^{l}(-1)^{j}{2l\choose l-j}S_{jt}f, \quad f\in L^2(S^m).
$$
Since the cosine function is defined in the whole real line, $S_{l,t}$ is well-defined while $S_{1,t}=S_t$.\
Taking advantage of the arguments delineated in the previous example and keeping the notation used there, one can see that
$$
S_{l,t} (f)=-2{2l\choose l}^{-1}\sum_{j=1}^{l}(-1)^{j}{2l\choose l-j}(f\ast \mu_{jt}^m), \quad f\in L^p(S^m).
$$
Consequently,
$$S_{l,t}(f)=f\ast \mu_t^m(l), \quad f \in L^2(S^m), \quad t\in (0,\pi),$$
where
$$
 \mu_t^m(l):= -2{2l\choose l}^{-1}\sum_{j=1}^{l}(-1)^{j}{2l\choose l-j}\mu_{jt}^m.
$$
The sequence of multipliers of $S_{l,t}$ is defined by
$$
\eta_k^t(l)=-2{2l\choose l}^{-1}\sum_{j=1}^{l}(-1)^{j}{2l\choose l-j}\frac{C_k^{(m-1)/2}(\cos jt)}{C_k^{(m-1)/2}(1)}, \quad l=1,2,\ldots, \quad k \in \mathbb{Z}_+, \quad t\in(0,\pi).
$$
Lemma 4.4 in \cite{dai} points that
$$
1-\eta_k^t(l)\asymp (\min\{1, kt\})^{2l}, \quad k,l=1,2,\ldots, \quad t\in(0,\pi/2),\quad  k \in \mathbb{Z}_+.
$$ Finally, the uniform boundedness of the family $\{S_{l,t}: t \in (0,\pi)\}$ follows from
the inequality
$$\|S_{l,t} f\|_2 \leq \left(2^{2l}{2l\choose l}^{-1}-1\right) \|f\|_2, \quad f \in L^2(S^m), \quad t \in (0,\pi).$$

\noindent  \textbf{Averages on caps.} This example has its origin in the reference \cite{berens} but is also discussed in details in \cite{ditzian2}.\ The average operator on the cap
$$
C_t^x=\{w\in S^m: x \cdot y\geq \cos t\}
$$
of $S^m$, defined by $t$ and the pole $x$, is the operator $M_t$ given by
$$
(A_t f)(x)=\frac{1}{C_m(t)}\int_{C_t^x}f(w)d\sigma_m(w), \quad x\in S^m,\quad  t\in (0,\pi),
$$
in which $dr$ corresponds to integration over $C_t^x$ and
$C_m(t)$ is total volume of the cap $C_t^x$.\ Since the right-hand side of
$$
C_m(t)=\omega_{m-1}\int_0^{t}(\sin h)^{m-1}dh, \quad t \in (0,\pi),
$$
does not depend upon $x$, the notation $C_m(t)$ is plainly justified.\ In order to see that $A_t$ is a convolution operator, we consider the auxiliary zonal kernel
$$\mathcal{Z}_{t}(x, y):=\left\{\begin{array}{rc}
\omega_m,& \quad \mbox{if\ \ } \cos t\leq x\cdot y\leq 1\\
0,& \quad otherwise.
\end{array}\right.$$
Clearly,
$$A_t(f)(x)=\frac{1}{\omega_m C_m(t)}\int_{S^m} \mathcal{Z}_t(x,y)f(y)d\sigma_m(y), \quad x \in S^m, \quad f \in L^2(S^m), \quad t \in (0,\pi),$$
that is, $A_t(f)$ is the usual spherical convolution of $C_m(t)^{-1}\mathcal{Z}_{t}$ with $f$ (\cite{berens}).\ To see that $A_t$ fits into the setting of spherical convolution with measures, it suffices to remember that
$L^2_\varepsilon(S^m)$ is embeddable in $\mathcal{M}_\varepsilon(S^m)$.\ The embedding itself is $f \mapsto \mu_f$ in which
$$
d\mu_f(x)=f(x)d\sigma_m(x), \quad x \in S^m.
$$
Since each $\mathcal{Z}_{t}(\varepsilon,\cdot)$
belongs to $L^2_\varepsilon(S^m)$, the embedding produces a corresponding element $\tilde{\mu}_{\mathcal{Z}_t}$ in $\mathcal{M}_\varepsilon(S^m)$ defined by
$$d\tilde{\mu}_{\mathcal{Z}_t}(x)=\mathcal{Z}_{t}(\varepsilon,x)d\sigma_m(x), \quad x \in S^m.$$ Finally, since $\mu_{\mathcal{Z}_{t}}=C_m^{-1}(t)\tilde{\mu}_{\mathcal{Z}_t}$ and
$$
f\ast \mu_{\mathcal{Z}_{t}}(x)=\frac{1}{\omega_m}\int_{S^m}f(y)d\varphi_x(\mu_{\mathcal{Z}_{t}})(y)=\frac{1}{\omega_m}\int_{S^m}f(y)d(\mu_{\mathcal{Z}_{t}} \circ \mathcal{O}_x^\varepsilon)(y),
$$
we immediately obtain
$$
f\ast \mu_{\mathcal{Z}_{t}}(x)=\frac{1}{\omega_m C_m(t)}\int_{S^m}f(y){\mathcal{Z}_{t}}(\varepsilon,\mathcal{O}_x^\varepsilon y)d\sigma_m(y)=\frac{1}{\omega_mC_m(t)}\int_{S^m}\mathcal{Z}_{t}(x,y)f(y)d\sigma_m(y),
$$
that is,
$$
M_t(f)=f\ast \mu_{\mathcal{Z}_{t}},\quad f \in L^2(S^m), \quad t \in (0,\pi).
$$
The sequence of multipliers of $M_t$ is (\cite{berens})
$$
\rho_k^{t}=\frac{\omega_{m-1}}{C_k^{(m-1)/2}(1)C_m(t)}\left(\int_{0}^{t}C_k^{(m-1)/2}(\cos h)(\sin h)^{m-1}dh\right), \quad  k\in \mathbb{Z}_+, \quad t\in(0,\pi).
$$
With the same constants in the shifting operator case, we have that
$$
\frac{\omega_{m-1}c_1k^2}{C_m(t)}\left(\int_{0}^{t}h^2(\sin h)^{m-1}dh\right)\leq 1-\rho_k^t\leq \frac{\omega_{m-1}c^2k^2}{C_m(t)}\left(\int_{0}^{t}h^2(\sin h)^{m-1}dh\right),\quad k \in \mathbb{Z}_+.
$$
After we estimate the sine function, this double inequality takes the form
$$
\frac{\omega_{m-1}c_1n^2t^{m+2}}{(2\pi)^{m-1}(m+2)C_m(t)} \leq 1-\rho_k^t \leq \frac{\omega_{m-1}c_2n^2t^{m+2}}{(m+2)C_m(t)}, \quad k \in \mathbb{Z}_+.
$$
On the other hand, direct computation yields
$$
\frac{\omega_{m-1}}{m}\left(\frac{2}{\pi}\right)^{m-1} t^m\leq C_m(t) \leq \omega_{m-1}t^m,  \quad t \in (0,\pi),
$$
so that
$$
c_1'(tk)^2\leq 1-\rho_k^t\leq c_2'(tk)^2 \quad 0 < kt \leq \pi, \quad t\in (0,\pi/2],
$$
for convenient positive constants $c_1'$ and $c_2'$.\ Similarly, still keeping the notation for the example involving the shifting operator, we can deduce that
for any $\tau>0$,
$$
0< \rho_k^t\leq c_{\tau}, \quad t\geq\tau/k,
$$
with $c_{\tau}<1$ depending on $m$ and $\tau$.\ Thus,
$$
1-\rho_k^t\asymp (\min\{1, kt\})^{2},  \quad k \in \mathbb{Z}_+, \quad t \in (0,\pi).
$$
Finally, the inequality (\cite{berens})
$$\|M_t f\|_2 \leq \|f\|_2, \quad f \in L^2(S^m), \quad t \in (0,\pi),$$
provides the uniform boundedness of the family $\{M_t: t \in (0,\pi)\}$.
\vspace*{2mm}

\noindent \textbf{Stekelov-type means.} This operator was introduced in \cite{ditzian2} and gives us an interesting additional example.\ The Stekelov-type mean is given by
$$
E_t(f)(x)=\frac{1}{D_m(t)}\int_0^{t}\frac{C_m(s)}{R_m(s)} A_s(f)(x)ds, \quad x\in S^m,\quad  t\in (0,\pi),
$$
the normalizing constant $D_m(t)$ being chosen so that $E_t(1)=1$.\ In order to see that the operators $E_t$ fit into the convolution structure we are using, let us consider the family of locally supported kernels defined by the formula
$$\mathcal{W}_{t}(x, y):=\left\{\begin{array}{rc}
\displaystyle{\int_0^t \frac{1}{R_m(s)}\mathcal{Z}_{s}(x,y)ds},& \quad \mbox{if\ \ } \cos t \leq x\cdot y\leq  1\\
0,& \quad otherwise
\end{array}\right.$$
where $\mathcal{Z}_{s}$ are the kernels described in the previous example.\ Since the kernel $\mathcal{W}_{t}$ is bizonal, the procedure developed before can be reproduced here in order to see that the formula
$$
d\tilde{\mu}_{\mathcal{W}_{t}}(x)=\mathcal{W}_{t}(\varepsilon\cdot x)d\sigma_m(x), \quad x \in S^m,
$$
defines an element $\tilde{\mu}_{\mathcal{W}_{t}}$ in $\mathcal{M}_\varepsilon(S^m)$ and, consequently, if $\mu_{\mathcal{W}_{t}}:=D_m^{-1}(t)\tilde{\mu}_{\mathcal{W}_{t}}$, then
$$
E_t(f)=f \ast \mu_{\mathcal{W}_{t}}, \quad f \in L^2(S^m),\quad t \in (0,\pi).
$$
Also, it is not difficult to see that
$$
1-\varphi_k^t\asymp (\min\{1, kt\})^{2}, \quad k \in \mathbb{Z}_+,\quad t \in (0,\pi),
$$
 where $\{\varphi_k^t\}$ is the multiplier family associated to $\{E_t: t\in (0,\pi)\}$ and given by
$$
\varphi_k^t=\frac{1}{D_m(t)}\int_0^{t}\frac{C_m(s)\rho_k^s}{R_m(s)}ds, \quad k \in \mathbb{Z}_+,\quad t\in (0,\pi).$$
The calculations in this case are similar to those done in the previous one, reason why we will not reproduce the details here.

\section{Remarks}

Most of the concepts and constructions made in this paper can be recovered when we replace the unit sphere with a compact symmetric space of rank 1.\ Indeed, these spaces are Riemannian manifolds possessing a harmonic analysis structure very similar to that we have on the spheres.\ A well-known classification for these spaces is as follows: the spheres $S^m$ ($m=1,2,\ldots$), the real projective spaces $P^m(\mathbb{R})$ ($m=2,3,\ldots$), the complex projective spaces $P^m(\mathbb{C})$ ($m=4,6,\ldots$), the quaternion projective spaces $P^m(\mathbb{H})$ ($m=8,12, 16\ldots$) and Cayley's elliptic plane $P^{16}$.\ Additional information about them and pertinent to a possible extension of the results in this paper can be found in papers authored by S. S. Platonov (for example, the survey paper \cite{platonov}).

The decay presented in Theorem \ref{mainn} and its corollary seems to be optimal within the setting considered.\ After some attempts, we were unable to find either an example or a decent argument substantiating such assertion.


%
%

\vspace*{2cm}

\noindent 
Departamento de
Matem\'atica,\\ ICMC-USP - S\~ao Carlos, Caixa Postal 668,\\
13560-970 S\~ao Carlos SP, Brasil\\ e-mails: thsjordao@gmail.com; menegatt@icmc.usp.br

\end{document}